\documentclass[12pt]{amsart}
\usepackage{amsmath,amsfonts,amstext,amsthm,epsfig}

\newtheorem{thm}{Theorem}[section]
\newtheorem{cor}[thm]{Corollary}
\newtheorem{lem}[thm]{Lemma}
\newtheorem{prop}[thm]{Proposition}
\newtheorem{defn}[thm]{Definition}
\theoremstyle{remark}

\numberwithin{equation}{section}

\newcommand{\weakto}{\rightharpoonup}

\def\cl{\overline}

\def\N{\mathbb N}
\def\Z{\mathbb Z}

\def\R{\mathbb R}

\begin{document}


\title[Conley index continuation under time averaging]{Averaging,
Conley index continuation\\
and recurrent dynamics in almost-periodic\\ parabolic equations}
\author[Martino Prizzi]{Martino Prizzi}
\address{Martino Prizzi, Dipartimento di Scienze Matematiche, Universit\`a degli
Studi di Trieste,
via Valerio 12/b, I-34127 Trieste, Italy}
\email{prizzi@mathsun1.univ.trieste.it}%


\begin{abstract}
We study a non-autonomous parabolic equation with
almost-periodic, rapidly oscillating principal part and nonlinear
interactions. We associate to the equation a skew-product semiflow
and, for a special class of nonlinearities, we define the Conley
index of an isolated  invariant set. As the frequency of the
oscillations tends to infinity, we prove that every isolated
invariant set of the averaged autonomous equation can be continued
to an isolated invariant set of the skew-product semiflow
associated to the non-autonomous equation. Finally, we illustrate
some examples in which the Conley index can be explicitely
computed and can be exploited to detect the existence of recurrent
dynamics in the equation.
\end{abstract}


\maketitle


\section{Introduction}
In this paper we study a family of non-autonomous parabolic
equations
\begin{equation} \label{intro1}
u_t-\sum_{i,j=1}^Na_{ij}(\omega
t)\partial_i\partial_j u=F(\omega
t,x,u),\quad(t,x)\in\R\times\R^N
\end{equation}
with almost-periodic, rapidly oscillating principal part and
nonlinear interactions.
Under suitable hypotheses (see Section 2), the Cauchy problem for
(\ref{intro1}) is well-posed in $H^1(\R^N)$ and the equation
generates a (local) process, that is a
two-parameter family of nonlinear operators $\Pi_\omega(t,s)$
such that $\Pi_\omega(t,t)=I$,
$t\in\R$, and $\Pi_\omega(t,p)\Pi_\omega(p,s)=\Pi_\omega(t,s)$, $t\geq p \geq s$.
\par We are interested in
the behaviour of the solutions of (\ref{intro1}) as $\omega
\rightarrow +\infty$. It is well known that, if a
function $\sigma$ is almost-periodic, then its mean
value
\begin{equation}\label{intro1.5}
\lim_{T\rightarrow +\infty}{{1}\over{2T}}\int_{-T}^T \sigma(p)\,dp =:\bar \sigma
\end{equation}
 is well defined.
This fact suggests that the averaged equation
\begin{equation}\label{intro2}
u_t-\sum_{i,j=1}^N \bar a_{ij}\partial_i\partial_j u= \bar F(x,u),
\quad (t, x)\in\R_+\times\R^N
\end{equation}
should behave like a limit equation for (\ref{intro1}) as
$\omega \to +\infty$.\par Results of this kind have been known for quite a long time for
ordinary differential equations with almost-periodic coefficients, and are related to the so called
{Bogolyubov averaging principle} (see \cite{Bogo}). For evolution equations
in infinite dymensions, {\it local} results in this direction have been
obtained in an abstract setting by Hale and Verduyn Lunel \cite{hallun}.
In a more recent paper \cite{il}, Ilyin
proposes a {\it global} criterion for comparison between the
process generated by an almost-periodic equation and the semiflow generated by
the corresponding averaged equation. The model problem is a parabolic equation
on a bounded domain, with
an almost-periodic time-dependent nonlinearity. Under suitable dissipativeness
and compactness hypotheses,
both the process and the semiflow possess compact global attractors (see \cite{chepvish}).
A first (rough) way to express the concept
of closeness of the two is then to give an estimate of the Hausdorff distance
of their attractors. A more detailed description of the internal
structure of the attractors is given by Efendiev and Zelik in \cite{efzel}.
They assume that the averaged problem admits a Lyapunov functional
and that the semiflow on the attractor is Morse-Smale. Then they show
that this structure, in a certain sense, persists in the almost-periodic
perturbation, provided the frequency of the oscillations is
sufficiently large.

\par The aim of this paper is to investigate the persistence, under almost-periodic and
rapidly oscillating perturbations,
of invariant sets which are possibly more general than attractors or hyperbolic equilibria.
This task leads naturally to the use of global topological tools like the homotopy index
of Conley.

\par Let $X$ be a metric space and let $\pi$ be a local semiflow in $X$. If $K$ is an isolated
$\pi$-invariant set for which there exists a
$\pi$-admissible isolating neighborhood $B$ (see \cite{ryba} for the precise
definitions of this and of the related concepts), then one can prove that there exists a
special isolating neighborhood
${\mathcal B}\subset B$ of
$K$, called an {\it isolating block}, which has the property that solutions of $\pi$ are
``transverse" to the boundary of
${\mathcal B}$. Letting ${\mathcal B}^-$ be the set of all points of $\partial{\mathcal B}$ the solutions
through which leave ${\mathcal B}$ in positive time direction, and collapsing ${\mathcal B}^-$ to one
point, we obtain the
{\it pointed space} ${\mathcal B}/{\mathcal B}^-$ with the distinguished {\it base point} $p=[{\mathcal B}^-]$.
It turns out that the homotopy type $h({\mathcal B}/{\mathcal B}^-,[{\mathcal B}^-])$
of $({\mathcal B}/{\mathcal B}^-,[{\mathcal B}^-])$ does not depend on the choice of ${\mathcal B}$.
This means that $h({\mathcal B}/{\mathcal B}^-,[{\mathcal B}^-])$ depends only on the pair $(\pi, K)$,
and we write $h(\pi,K):=h({\mathcal B}/{\mathcal B}^-,[{\mathcal B}^-])$. $h(\pi,K)$ is called the
{\it homotopy index} of $(\pi,K)$. For two-sided
flows on locally compact spaces, the homotopy index is due to Charles Conley
(see  \cite{conley}) and therefore it is called the {\it Conley index}. In the case of
a local semiflow $\pi$ in an arbitrary metric space $X$, the extended homotopy index theory
was developed by Rybakowski in \cite{ryba1} and rests in an essential way on the
notion of {\it $\pi$-admissibility}. The most important properties of the Conley index
are the following: (a) if
$h(\pi,K)\not=\underline 0$, then
$K\not=\emptyset$; (b) the homotopy index is invariant under continuation, in the
sense that, roughly speaking, it remains constant along ``continuous"  deformations of the
pair $(\pi,K)$.

\par The first difficulty in applying the homotopy index theory to (\ref{intro1}) comes
from the fact that non-autonomous equations define {\it processes} and not {\it semiflows}.
The theory of {\it skew-product semiflows}, developed by Sell in
\cite{sell}, provides then the right functional setting for a dynamical-system treatment of
equation (\ref{intro1}), at the expense of introducing an extended phase space. Another
difficulty
comes from the characteristic {\it lack of compactness} exibited by problems in {\it
unbounded domains}. In fact, in the case of a parabolic
equation on a {\it bounded} open set $\Omega\subset\R^N$,
the admissibility of
all bounded closed sets in the phase space is a direct consequence of the compactness of the
Sobolev embedding $H^1(\Omega)\hookrightarrow L^2(\Omega) $. In $\R^N$ this property fails,
and one has to introduce some restrictions on the non-linear term
$F$. The question of admissibility for autonomous equations in unbounded domains was discussed in
\cite{priz}, where a condition on $F$ was given, ensuring
the admissibility of all bounded closed sets in the phase space.
In the same spirit, we shall assume here that the nonlinearity $F$
satisfies a condition like
\begin{equation}\label{intro3}
F(\tau,x,u)u\leq-\nu
|u|^2+b(\tau,x)|u|^q+c(\tau,x),
\end{equation}
where $b(\tau,x)$ and
$c(\tau,x)$ tend to $0$ as $|x|\to\infty$, in some sense to be
made precise later. Roughly speaking, (\ref{intro3}) means that the
nonlinearity $F$ is dissipative for large $x$. Therefore, we term (\ref{intro3}) as a
``dissipativeness-in-the-large" condition.

\par It seems that the first to use the homotopy index in connection with the averaging priciple
was Ward in \cite{ward1}. He considered an {\it ordinary
differential equation} with non-autonomous, almost-periodic
nonlinearity. He proved that if the autonomous averaged equation
possesses an isolated invariant set with nontrivial homotopy
index, the latter can be continued to  a nearby isolated invariant
set of the skew-product flow associated to the non-autonomous
equation, provided the frequency of the oscillations is
sufficiently large. From this he deduced the existence of {\it
bounded full solutions} of the original non-autonomous equation.

\par In this paper we procede in a similar way. We define a skew-product semiflow in the space
$\Sigma\times H^1(\R^N)$, where $\Sigma$ is the ``symbol space" associated to the non-autonomous
equation (\ref{intro1}).   Then we prove that, under the ``dissipativeness-in-the-large" condition
(\ref{intro3}), all bounded closed sets in the extended phase space are admissible. Therefore it is possible
to define the Conley index of an isolated  invariant set. As the frequency of the oscillations tends to
infinity, we prove that every isolated invariant set of the averaged autonomous equation can be
continued to an isolated invariant set of the skew-product semiflow associated to the non-autonomous
equation. Again, from this we can easily deduce  the existence of {\it bounded full solutions} of the
original non-autonomous equation. However, from the dynamical point of view, it is
much more interesting to look for {\it recurrent solutions} (in the sense of Birkhoff) rather than for
{\it bounded solutions} of the equation (\ref{intro1}). In the last section, we briefly recall the
concept of {\it recurrence} and we show that, under a technical condition on the principal coefficients
$a_{ij}(\cdot)$, the existence of {\it recurrent} solutions of (\ref{intro1}) is a straightforward
consequence of the existence of a non-empty, compact invariant set of the corresponding skew-product
semiflow. We conclude with an example, in which the averaged equation is asymptotically linear and
the homotopy index can be explicitly computed.

\section{The process and its properties}

We consider the non-autonomous parabolic equation
\begin{equation}\label{equazione}
u_t-\sum_{i,j=1}^N a_{ij}(\omega t)\partial_i\partial_j u=F(\omega
t,x,u),\end{equation}
where $(t,x)\in\R\times\R^N$ and
$\omega$ is a positive constant. \par For notational convenience,
we shall assume throughout that $N\geq3$. We make the following
assumptions:

\begin{itemize}
\item[\bf(H1)] {\it for every $\tau\in\R$ the matrix $(a_{ij}(\tau))_{ij}$ is real
symmetric. There exists a constant $\nu_0>0$ such that $\nu_0|\xi|^2\leq
\sum_{ij}a_{ij}(\tau)\xi_i\xi_j\leq \nu_0^{-1}|\xi|^2$ for all $(\tau,\xi)\in\R\times\R^N$.
There exist a constant $0<\theta<1$ and a positive constant $C$ such that, for all
$\tau_1,\tau_2\in\R$, and for $1\leq i,j\leq N$,
\begin{equation}\label{a0}
|a_{ij}(\tau_1)-a_{ij}(\tau_2)|\leq C|\tau_1-\tau_2|^\theta;
\end{equation}
}
\item[\bf(H2)] {\it the function $F$ is continuous on $\R\times\R^N\times\R$ and
for every
$\tau\in\R$ the function
$F(\tau,\cdot,0)$ is square integrable;}

\item[\bf(H3)] {\it for every $(\tau,x)\in\R\times\R^N$ the function
$F(\tau,x,\cdot)$ is continuously differentiable
 and there exists a constant $C$ such that
\begin{equation}\label{a1}
|F'_u(\tau,x,u)|\leq C(1+|u|^\beta)\quad\text{for all
$(\tau,x,u)\in\R\times\R^N\times\R$,}\end{equation}
where $\beta:= 2^*/2-1$;}

\item[\bf(H4)] {\it there exist a constant $0<\theta<1$, a positive constant $C$ and a
function
$g_0\in L^2(\R^N)$ such that, for all $\tau_1,\tau_2\in\R$ and $(x,u)\in\R^N\times\R$,
\begin{equation}\label{a2}
|F(\tau_1,x,u)-F(\tau_2,x,u)|\leq C(g_0(x)+|u|+|u|^{\beta+1})|\tau_1-\tau_2|^\theta.
\end{equation}
}
\end{itemize}

Let ${\mathcal M}_1$ be the space of $N\times N$ real symmetric matrices and define ${\mathcal M}_2$
to be the space of all functions
$f\colon\R^N\times\R\to\R$ such that
$f(x,u)$ satisfies {\bf(H1)} and {\bf(H2)},  equipped with
the norm
\begin{equation}\label{a3}
\|f\|_{{\mathcal M}_2}
:=\|f(\cdot,0)\|_{L^2}+\sup_{(x,u)\in\R^N\times\R}(1+|u|^\beta)^{-1}|f'_u(x,u)|.
\end{equation}
We assume that
\begin{itemize}
\item[\bf(AP)] {\it the functions $\tau\mapsto (a_{ij}(\tau))_{ij}\in{\mathcal M}_1$ and
$\tau\mapsto
F(\tau,\cdot,\cdot)\in{\mathcal M}_2$ are almost-periodic.}
\end{itemize}

We recall some basic facts on almost-periodic functions. By Bochner's criterion (see e.g.
\cite{levzhik}), whenever ${\mathcal M}$ is a Banach space and
$\sigma\colon \R\to{\mathcal M}$ is almost-periodic, the set of all translations
$\{\,\sigma(\cdot+h)\mid h\in\R\,\}$ is precompact in
$C_b(\R,{\mathcal M})$. The closure of this set in $C_b(\R,{\mathcal M})$ is called the hull of
$\sigma$
and is usually denoted by ${\mathcal H}(\sigma)$. Moreover, if $\zeta\in {\mathcal H}(\sigma)$, then
$\zeta$ is almost-periodic and ${\mathcal H}(\zeta)={\mathcal H}(\sigma)$. We recall also
that, for an almost-periodic function $\sigma$, the mean value
\begin{equation}
\lim_{T\to\infty}{1\over{2T}}\int_{-T}^{T}\sigma(t)\, dt=\bar\sigma\in{\mathcal M}
\end{equation}
exists. More remarkably, one can prove (see again \cite{levzhik}) that there exists a
bounded decreasing function
$\mu\colon \R_+\to\R_+$, $\mu(T)\to 0$ as $T\to\infty$, such that
\begin{equation}\label{media}
\|(1/T)\int_s^{s+T}(\zeta(t)-\bar\sigma)\,dt\|_{{\mathcal M}}\leq\mu(T)\quad
\text{for all $s\in\R$ and all $\zeta \in {\mathcal H}(\sigma)$.}
\end{equation}
If ${\mathcal M}$, ${\mathcal N}$ are Banach spaces and
$\sigma\colon \R\to{\mathcal M}$, $\rho\colon \R\to{\mathcal N}$ are almost-periodic,
then $(\sigma,\rho)\colon\R\to{\mathcal M}\times{\mathcal N}$ is almost-periodic
and ${\mathcal H}((\sigma,\rho))\subset{\mathcal H}(\sigma)\times{\mathcal H}(\rho)$.
Moreover, the mean value of
$(\sigma,\rho)$ is $(\bar\sigma,\bar\rho)$.

We denote by $\Sigma_1$ and $\Sigma_2$ the hulls
of the functions $\tau\mapsto (a_{ij}(\tau))_{ij}$ and  $\tau\mapsto F(\tau,\cdot,\cdot)$
in $C_b(\R,{\mathcal M}_1)$ and
$C_b(\R,{\mathcal M}_2)$ respectively.
The corresponding mean values are denoted by $(\bar a_{ij})_{ij}\in{\mathcal M}_1$ and
$\bar F(\cdot,\cdot)\in {\mathcal M}_2$. Besides, we denote by $\Sigma$ the hull of
$\tau\mapsto ((a_{ij}(\tau))_{ij}, F(\tau,\cdot,\cdot))$ in $C_b(\R,{\mathcal M}_1\times{\mathcal
M}_2)$. Sometimes $\Sigma$ is called the ``symbol space" associated to the equation.

It is easy
to check that {\bf(H1)} is satisfied
by any element of $\Sigma_1$ as well as by
the corresponding mean value, and {\bf(H2)}--{\bf(H4)} are satisfied
by any element of $\Sigma_2$ as well as by
the corresponding mean value (with the same
constants).

For later use, we need also to introduce a parameter $\lambda\in[0,1]$.
For $\lambda\in[0,1]$ and $((\alpha_{ij}(\cdot))_{ij},\Phi(\cdot,\cdot,\cdot))\in\Sigma$,
we define
\begin{equation}
\alpha_{ij}(\lambda,\tau):=\lambda\alpha_{ij}(\tau)+(1-\lambda)\bar a_{ij},\quad 1\leq
i,j\leq N
\end{equation}
and
\begin{equation}
\Phi(\lambda,\tau,x,u):=\lambda\Phi(\tau,x,u)+(1-\lambda)\bar F(x,u),\quad
(\tau,x,u)\in\R\times\R^N\times\R.
\end{equation}
Notice that $\tau\mapsto\alpha_{ij}(\lambda,\tau)$ and
$\tau\mapsto\Phi(\lambda,\tau,\cdot,\cdot)$ are  almost-periodic and their mean values are
respectively $(\bar a_{ij})_{ij}$ and $\bar F(\cdot,\cdot)$.

 We introduce the Nemitski
operator
\begin{equation*}
\hat\Phi(\lambda,\cdot,\cdot)\colon\R\times H^1(\R^N)\to L^2(\R^N)
\end{equation*}
defined by
\begin{equation}
\hat\Phi(\lambda,\tau,u)(x):=\Phi(\lambda,\tau,u(x)).
\end{equation}
The map $\hat\Phi$ is continuos on $[0,1]\times\R\times H^1(\R^N)$ and differentiable with
respect to $u\in H^1(\R^N)$, and the following estimates hold:
\begin{equation}\label{est1}
\|\hat\Phi(\lambda,\tau,u)\|_{L^2}\leq C(1+\|u\|_{H^1}^{\beta+1}),
\end{equation}
\begin{equation}\label{est2}
\|D\hat\Phi(\lambda,\tau,u)\|_{{\mathcal L}(L^2,H^1)}\leq C(1+\|u\|_{H^1}^{\beta})
\end{equation}
and
\begin{multline}\label{est3}
\|\hat\Phi(\lambda,\tau_1,u_1)-\hat\Phi(\lambda,\tau_2,u_2)\|_{L^2}\leq
+C(1+\|u_1\|_{H^1}^{\beta+1}+\|u_2\|_{H^1}^{\beta+1})|\tau_1-\tau_2|^\theta\\
+C(1+\|u_1\|_{H^1}^{\beta}+\|u_2\|_{H^1}^{\beta})\|u_1-u_2\|_{H^1},
\end{multline}
where $C$ is a positive constant, $\beta$ is the exponent of {\bf(H2)} and $\theta$ is the
H\"older exponent of {\bf(H4)}.

For $t\in\R$, $\lambda\in[0,1]$, $\alpha=(\alpha_{ij}(\cdot))_{ij}\in\Sigma_1$ and
$\omega>0$, we define the operator
$A^\alpha_{\lambda,\omega}(t)\colon H^2(\R^N)\to L^2(\R^N)$ by
\begin{equation}
A^\alpha_{\lambda,\omega}(t)u:=-\sum_{i,j=1}^N \alpha_{ij}(\lambda,\omega
t)\partial_i\partial_j u,\quad u\in H^2(\R^N).
\end{equation}
Then $A^\alpha_{\lambda,\omega}(t)$ is a self-adjoint positive operator in $L^2(\R^N)$ and
our assumptions on the coefficients $a_{ij}(\tau)$ imply that the abstract parabolic
equation
\begin{equation}
\dot u=-A^\alpha_{\lambda,\omega}(t)u
\end{equation}
generates a linear process
\begin{equation*}
U^\alpha_{\lambda,\omega}(t,s)\colon L^2(\R^N)\to L^2(\R^N), \quad t\geq s,
\end{equation*}
such that
\begin{equation}\label{expest}
\|U^\alpha_{\lambda,\omega}(t,s)u\|_{L^2}\leq M\|u\|_{L^2},\quad u\in L^2(\R^N),
\end{equation}
\begin{equation}\label{expest2}
\|U^\alpha_{\lambda,\omega}(t,s)u\|_{H^1}\leq M\|u\|_{H^1},\quad u\in H^1(\R^N),
\end{equation}
and
\begin{equation}\label{expest3}
\|U^\alpha_{\lambda,\omega}(t,s)u\|_{H^1}\leq M(1+(t-s)^{-1/2})\|u\|_{L^2},\quad u\in
L^2(\R^N),
\end{equation}
where $M$ is a positive constant depending only on $\nu_0$
(see e.g. \cite{Pazy}, Ch.5, and \cite{tana}).

\par
For $\lambda=0$, $\alpha_{ij}(\lambda,\omega t)\equiv\bar a_{ij}$. We set
$\bar A:=A^\alpha_{0,\omega}(t)$, so we have $U^\alpha_{0,\omega}(t,s)\equiv e^{-\bar
A(t-s)}$. Representing
$U^\alpha_{\lambda,\omega}(t,s)$ in terms of its Fourier transform,
one can prove (cf \cite{antopriz}, Propositions 4.1 -- 4.3) that $U^{\alpha}_{\lambda,\omega}(t,s)$
converges to $e^{-\bar A(t-s)}$ in a strong sense, uniformly with respect to $\alpha$ and $\lambda$.

For every $\lambda\in[0,1]$ and
$\sigma:=((\alpha_{ij}(\cdot))_{ij},\Phi(\cdot,\cdot,\cdot))\in\Sigma$, one can consider
the nonlinear equation
(\ref{equazione}) with $a_{ij}(\omega t)$ and $F(\omega t,x,u)$ replaced by
$\alpha_{ij}(\lambda,\omega t)$ and $\Phi(\lambda,\omega t,x,u)$ respectively. Following
\cite{He}, we rewrite equation
(\ref{equazione})  as an abstract evolution equation, namely
\begin{equation}\label{pbcauchy}
\begin{cases}
\dot u+A^\alpha_{\lambda,\omega}(t)u=\hat \Phi(\lambda,\omega t,u)\\
u(s)=u_s
\end{cases}
\end{equation}

By classical results of \cite{fried}, \cite{He} and \cite{Pazy},
for every
$s\in\R$ and  $u_s\in H^1(\R^N)$, the semilinear Cauchy
problem (\ref{pbcauchy})
is locally well-posed. More  specifically, one has the following

\begin{prop} \label{localexist} For every
$R>0$ there exists
$T_R>0$ (independent of $s$, $\sigma$, $\lambda$ and $\omega$)
such that, for all $u_s\in B_{H^1}(R;0)$, problem (\ref{pbcauchy}) admits a unique solution
$u(\cdot)$ defined for
$t\in[s,s+T_R]$, with $\|u(t)\|_{H^1}\in B_{H^1}(2R;0)$.\end{prop}

It follows that problem (\ref{pbcauchy}) possesses a unique
maximal solution
$u\in C^0([s,s+T[, H^1)\cap C^1(]s,s+T[,L^2)$, where $T$ depends on  $u_s$.
The solution
$u(\cdot)$ satisfies the variation-of-constant formula
\begin{equation}\label{varconst1}
u(t)=U^\alpha_{\lambda,\omega}(t,s)u_s+\int_s^t U^\alpha_{\lambda,\omega}(t,p)
\hat\Phi(\lambda,\omega p,u(p))\,dp,\quad t\geq s.
\end{equation}

It follows that for every $\lambda\in[0,1]$ and
$\sigma:=((\alpha_{ij}(\cdot))_{ij},\Phi(\cdot,\cdot,\cdot))\in\Sigma$, equation (\ref{pbcauchy})
generates a local process $\Pi^\sigma_{\lambda,\omega}(t,s)$.

Thanks to the variation-of-constant formula (\ref{varconst1}), one can prove
(cf \cite{antopriz}, Lemma 3.6) the following

\begin{lem} \label{contin1}
Let $\sigma\in\Sigma$ and let $(\sigma_n)_{n\in\N}$ be a sequence in
$\Sigma$, such that $\sigma_n\to\sigma$ as $n\to\infty$.
Let $\lambda\in[0,1]$ and let $(\lambda_n)_{n\in\N}$ be a sequence in
$[0,1]$, such that $\lambda_n\to\lambda$ as $n\to\infty$.
Let $u\in H^1(\R^N)$ and let
$(u_n)_{n\in\N}$ be a  bounded sequence in $H^1(\R^N)$. Let $T>0$ and
let $(t_n)_{n\in\N}$ and
$(s_n)_{n\in\N}$ be two sequences of real numbers, with $t_n\in[s_n,s_n+T]$ for all $n$, and
assume that $t_n\to t$ and $s_n\to s$ as $n\to\infty$. Let $\omega>0$. Finally, let $R>0$ and
assume that, for all $n$, $\|\Pi^{\sigma_n}_{\lambda_n,\omega}(r,s_n)u_n\|_{H^1}\leq R$,
$r\in[s_n,s_n+T]$,  and $\|\Pi^{\sigma}_{\lambda,\omega}(r,s)u\|_{H^1}\leq R$,
$r\in[s,s+T]$.  Then
\begin{enumerate}
\item if $u_n\to u$ in $L^2(\R^N)$ and $t>s$,
$$
\|\Pi^{\sigma_n}_{\lambda_n,\omega}(t_n,s_n)u_n-\Pi^\sigma_{\lambda,\omega}(t,s)u\|_{H^1}\to
0\quad\text{as
$n\to\infty$;}
$$
\item if $u_n\to u$ in $H^1(\R^N)$ and $t\geq s$,
$$
\|\Pi^{\sigma_n}_{\lambda_n,\omega}(t_n,s_n)u_n-\Pi^\sigma_{\lambda,\omega}(t,s)u\|_{H^1}\to
0\quad\text{as
$n\to\infty$.}
$$
\end{enumerate}
\end{lem}

A direct consequence of the second part of Lemma \ref{contin1} is
the following

\begin{prop}\label{1contin1}
Let $\sigma\in\Sigma$ and let $(\sigma_n)_{n\in\N}$ be a sequence in
$\Sigma$, such that $\sigma_n\to\sigma$ as $n\to\infty$.
Let $\lambda\in[0,1]$ and let $(\lambda_n)_{n\in\N}$ be a sequence in
$[0,1]$, such that $\lambda_n\to\lambda$ as $n\to\infty$.
Let $u\in H^1(\R^N)$ and let
$(u_n)_{n\in\N}$ be a  bounded sequence in $H^1(\R^N)$, such that $u_n\to u$ in $H^1(\R^N)$
as $n\to\infty$.
Let $(t_n)_{n\in\N}$ and
$(s_n)_{n\in\N}$ be two sequences of real numbers, and
assume that $t_n\to t$ and $s_n\to s$ as $n\to\infty$. Let $\omega>0$. Finally,
assume that $\Pi^{\sigma}_{\lambda,\omega}(r,s)u$ is defined for $r\in[s,t]$. Then, for all
$n$ sufficiently large,
$\Pi^{\sigma_n}_{\lambda_n,\omega}(r,s_n)u_n$ is defined for
$r\in[s_n,t_n]$  and
$$
\|\Pi^{\sigma_n}_{\lambda_n,\omega}(t_n,s_n)u_n-\Pi^\sigma_{\lambda,\omega}(t,s)u\|_{H^1}\to
0\quad\text{as
$n\to\infty$.}
$$
\end{prop}

For $\lambda=0$, (\ref{pbcauchy}) reduces to the autonomous problem
\begin{equation}\label{pbcauchyaut}
\begin{cases}
\dot u+\bar A u=\Hat{\Bar F}(u)\\
u(0)=u_0
\end{cases}
\end{equation}
For every
$u_0\in H^1(\R^N)$, the semilinear Cauchy
problem (\ref{pbcauchyaut})
is locally well-posed and hence possesses a unique maximal solution $u\in
C^0([0,T[, H^1)\cap C^1(]0,T[,L^2)$, where $T$ depends on $u_0$. Moreover,  $u$
satisfies the variation-of-constant formula
\begin{equation}\label{varconst2}
u(t)=e^{-\bar A t}u_0+\int_0^t e^{-\bar A(t-p)}
\Hat{\Bar F}(u(p))\,dp,\quad t\geq 0.
\end{equation}
The Cauchy problem (\ref{pbcauchyaut}) generates a local semiflow $\pi(t)$, and we have
$\Pi^\sigma_{0,\omega}(t,s)\equiv\pi(t-s)$.

By slightly modifying the proof of Theorem 4.4 in \cite{antopriz}, one can
prove the following averaging principle:

\begin{thm}\label{contin2}
Let $(\sigma_n)_{n\in\N}$ be a sequence in
$\Sigma$.
Let $(\lambda_n)_{n\in\N}$ be a sequence in
$[0,1]$.
Let $u\in H^1(\R^N)$ and let
$(u_n)_{n\in\N}$ be a  bounded sequence in $H^1(\R^N)$. Let $T>0$ and
let $(t_n)_{n\in\N}$ and
$(s_n)_{n\in\N}$ be two sequences of real numbers, with $t_n\in[s_n,s_n+T]$ for all $n$, and
assume that $t_n\to t$ and $s_n\to s$ as $n\to\infty$. Let $(\omega_n)_{n\in\N}$ be a
sequence of positive numbers, $\omega_n\to+\infty$ as
$n\to\infty$. Finally, let $R>0$ and
assume that, for all $n$, $\|\Pi^{\sigma_n}_{\lambda_n,\omega_n}(r,s_n)u_n\|_{H^1}\leq R$,
$r\in[s_n,s_n+T]$,  and $\|\pi(r-s)u\|_{H^1}\leq R$,
$r\in[s,s+T]$.  Then
\begin{enumerate}
\item if $u_n\to u$ in $L^2(\R^N)$ and $t>s$,
$$
\|\Pi^{\sigma_n}_{\lambda_n,\omega_n}(t_n,s_n)u_n-\pi(t-s)u\|_{H^1}\to
0\quad\text{as
$n\to\infty$;}
$$
\item if $u_n\to u$ in $H^1(\R^N)$ and $t\geq s$,
$$
\|\Pi^{\sigma_n}_{\lambda_n,\omega_n}(t_n,s_n)u_n-\pi(t-s)u\|_{H^1}\to
0\quad\text{as
$n\to\infty$.}
$$
\end{enumerate}
\end{thm}

Following \cite{chepvish}, we introduce the extended
phase-space $\Sigma\times H^1(\R^N)$. For $\omega>0$, we define on
$\Sigma$ the unitary group of translations
\begin{equation}\label{shift}
(T_\omega(h)\sigma)(\cdot):=\sigma(\cdot+\omega h).
\end{equation}
One can easily prove the following translation
identity:
\begin{equation}\label{transliden}
\Pi^\sigma_{\lambda,\omega}(t+h,s+h)=\Pi^{T_\omega(h)\sigma}_{\lambda,\omega}(t,s),
\quad h\in\R.
\end{equation}
Thanks to (\ref{transliden}), we can associate to the family of processes
$\{\,\Pi^\sigma_{\lambda,\omega}\mid\sigma\in\Sigma\,\}$ a
skew-product semiflow $P_{\lambda,\omega}(t)$ on the extended
phase-space $\Sigma\times H^1(\R^N)$, by the formula
\begin{equation}
P_{\lambda,\omega}(t)(\sigma,u):=(T_\omega(t)\sigma,\Pi_{\lambda,\omega}^\sigma(t,0)u).
\end{equation}
If $\omega>0$ and $\lambda\in[0,1]$ are fixed, Proposition \ref{localexist}
 implies  that the semiflow $P_{\lambda,\omega}$
satisfies the no-blow-up condition I-2.1 of \cite{ryba}.
Moreover, if $\omega>0$ is fixed and $(\lambda_n)_{n\in\N}$ is a
sequence converging to some $\lambda\in[0,1]$, Proposition \ref{1contin1}
 implies that the sequence of semiflows
$(P_{\lambda_n,\omega})_{n\in\N}$ converges to the semiflow
$P_{\lambda,\omega}$ on $\Sigma\times H^1(\R^N)$, according to
Definition I-2.2 of \cite{ryba}.
 Notice that, for
$\lambda=0$, one has
$P_{0,\omega}(t)(\sigma,u)=(T_\omega(t)\sigma,\pi(t)u)$, so
$P_{0,\omega}(t)(\sigma,u)$ is completely decoupled.

\section{The question of admissibility}
We begin by recalling the following concept, introduced by Rybakowski in \cite{ryba1}
(see also \cite{ryba}):
\begin{defn}\label{amm}
Let $X$ a metric space, let $B$ be a closed subset of $X$ and let $(\pi_n)_{n\in\N}$ be a
sequence of local semiflows in $X$. Then $B$ is called {\rm $\{\pi_n\}$-admissible} if the
following holds:
\par if $(x_n)_{n\in\N}$ is a sequence in $X$ and $(t_n)_{n\in\N}$ is a sequence in $\R_+$
such that $t_n\to\infty$ as $n\to\infty$ and $\pi_n(r)x_n\subset B$ for $r\in[0,t_n]$ for
all
$n\in\N$, then the sequence of endpoints $(\pi_n (t_n)x_n)_{n\in\N}$ has a converging
subsequence.
\par The set $B$ is called {\rm strongly $\{\pi_n\}$-admissible} if $B$ is $\{\pi_n\}$-admissible
and if $\pi_n$ does not explode in $B$ for every $n\in\N$. If $\pi_n=\pi$ for all $n$,
we say
that $B$ is {\rm $\pi$-admissible} (resp. {\rm strongly $\pi$-admissible})
\end{defn}

Notice that, by Proposition \ref{localexist}, if $B\subset H^1(\R^N)$
is bounded,
then the semiflow $P_{\lambda,\omega}$ does not explode in $\Sigma\times B$.
\par In the case of a parabolic
equation on a {\it bounded} open set $\Omega\subset\R^N$,
the admissibility of
all bounded subsets in the phase space is a direct consequence of the compactness of the
Sobolev embedding $H^1(\Omega)\hookrightarrow L^2(\Omega) $. In $\R^N$ this property fails,
and one has to introduce some restrictions on the non-linear term  $F$.
We make the following ``dissipativeness in the large" assumption (cf \cite{priz}):
\begin{itemize}
\item[\bf(D)] {\it for every $(\tau,x,u)\in\R\times\R^N\times\R$,
\begin{equation}\label{a4}
F(\tau,x,u)u\leq-\nu |u|^2+b(\tau,x)|u|^q+c(\tau,x),
\end{equation}
where $\nu>0$, $2\leq q<2N/(N-2)$, and $\tau\mapsto c(\tau,\cdot)\in L^1(\R^N)$ and
$\tau\mapsto b(\tau,\cdot)\in L^p(\R^N)$ are almost-periodic, where $2N/[2N-q(N-2)]\leq
p<\infty$.}
\end{itemize}
It is easy
to check that {\bf(D)} is satisfied
by any element of $\Sigma_2$ (with $b(\cdot,\cdot)$ and
$c(\cdot,\cdot)$ replaced by suitable functions $\beta(\cdot,\cdot)$ and
$\gamma(\cdot,\cdot)$ belonging to
the corresponding hulls) as well as by the mean value $\bar F$ (with $b(\cdot,\cdot)$ and
$c(\cdot,\cdot)$  replaced by their means $\bar b(\cdot)$ and $\bar c(\cdot)$).
Since the range of an almost-periodic function is compact, there exists a sequence  of
positive numbers
$(m_k)_{k\in\N}$,
$m_k\to0$ as
$k\to\infty$, such that
\begin{equation}\label{unifdecay}
\int_{|x|\geq k}|\beta(\tau,x)|^p\,dx+\int_{|x|\geq k}|\gamma(\tau,x)|\,dx\leq
m_k,\quad\tau\in\R,k\in\N,
\end{equation}
for all $\beta(\cdot,\cdot)\in{\mathcal H}(b(\cdot,\cdot))$ and
$\gamma(\cdot,\cdot)\in{\mathcal H}(c(\cdot,\cdot))$. Moreover,
\begin{equation}\label{unifdecay2}
\int_{|x|\geq k}|\bar b(x)|^p\,dx+\int_{|x|\geq k}|\bar c(x)|\,dx\leq
m_k,\quad k\in\N.
\end{equation}

The following Proposition is a non-autonomous version of Proposition 2.2 in \cite{priz}, and like the latter,
it was inspired by Lemma 5 in \cite{wang}:

\begin{prop}\label{comp1} Assume $(a_{ij}(\tau))_{ij}$ satisfies condition
{\bf(H1)}  and
$F(\tau,x,u)$ satisfies conditions {\bf(H2)}--{\bf(H4)}, {\bf(AP)} and {\bf(D)}. Let $R>0$.
There exists a sequence $(\eta_k)_{k\in\N}$, $\eta_k\to 0$ as
$k\to\infty$, with the following property:
\par whenever
$\lambda\in[0,1]$,
$\omega>0$, $(\alpha,\Phi)\in\Sigma$ and $u\colon [s,s+T]\to H^1(\R^n)$ is a solution
of (\ref{pbcauchy}) with
$\|u(t)\|_{H^1}\leq R$ for
$t\in[s,s+T]$, then
\begin{equation}
\int_{|x|\geq k}|u(t,x)|^2\,d x\leq R^2 e^{-2\nu (t-s)}+\eta_k\quad \text{ for
$t\in[s,s+T]$ and $k\in\N$.}
\end{equation}
The number
$\eta_k$ depends only on $R$, $\nu$, $\nu_0$ and $m_k$.
\end{prop}

\begin{proof}
Let $\theta\colon \R_+ \to \R$ be a smooth function such that $0\leq \theta(s) \leq 1$ for $s
\in \R_+$, $\theta(s)=0$ for $0\leq s \leq 1$ and $\theta(s)=1$ for $s \geq 2$.
Let $D := \sup_{s \in \R_+} |\theta'(s)|$. Define $\theta_k(x):=\theta(|x|^2/k^2)$. Then,
for $t\in[s,s+T]$, we have
\begin{multline*}
{{d}\over{dt}}{1\over2}\int_{\R^n}
\theta_k(x)|u(t,x)|^2\, d x =
\int_{\R^n}
\theta_k(x)u(t,x)u_t(t,x)\, d x \\
=-\int_{\R^n}\sum_{i,j=1}^N\alpha_{ij}(\lambda,\omega t)
\partial_i(\theta_k(x)u(t,x))\partial_j u(t,x)\, d x\\
+\int_{\R^n}\theta_k(x)u(t,x)\Phi(\lambda,\omega t,x,u(t,x))\,d x
\end{multline*}
Now we have
\begin{multline*}
-\int_{\R^n}\sum_{i,j=1}^N\alpha_{ij}(\lambda,\omega t)
\partial_i(\theta_k(x)u(t,x))\partial_j u(t,x)\, d x\\
=-\int_{\R^n}
\theta_k(x)\,\sum_{i,j=1}^N\alpha_{ij}(\lambda,\omega t)
\partial_iu(t,x)\partial_j u(t,x) d x\\
-{2\over k^2}\int_{\R^n}
\theta'(|x|^2/k^2)u(t,x)\sum_{i,j=1}^N\alpha_{ij}(\lambda,\omega t)\,x_i\,\partial_j u(t,x)\,
d x\\
\leq {2D\over {\nu_0k^2}}\int_{k\leq|x|\leq \sqrt 2 k}
|x|\, |u(t,x)|\,|\nabla_x u(t,x)|\, d x\leq {2\sqrt2D\over {\nu_0k}}R^2.
\end{multline*}
On the other hand, by condition {\bf(D)}, by the Sobolev embedding $H^1\hookrightarrow
L^{2n/(n-2)}$ and by H\"older inequality, we have
\begin{multline*}
\int_{\R^n}\theta_k(x)u(t,x)\Phi(\lambda,\omega t,x,u(t,x))\,d x
\leq -\nu \int_{\R^n}\theta_k(x)|u(t,x)|^2\,d x\\
+\int_{\R^n}\theta_k(x)(\lambda\beta(\omega t,x)+(1-\lambda)\bar
b(x))|u(t,x)|^q\,dx\\
+\int_{\R^n}\theta_k(x)(\lambda\gamma(\omega t,x)+(1-\lambda)\bar c(x))\,d x\\
\leq  -\nu \int_{\R^n}\theta_k(x)|u(t,x)|^2d x
+\left[{{(n-1)R}\over{(n-2)/2}}\right]^qm_k^{1/ p} +
m_k.
\end{multline*}
Summing up, we have found a sequence $(\eta_k)_{k\in\N}$, $\eta_k\to 0$ as $k\to\infty$,
such that
$$
{{d}\over{dt}}\int_{\R^n}
\theta_k(x)|u(t,x)|^2\, d x\leq -2\nu \int_{\R^n}\theta_k(x)|u(t,x)|^2\,d x+\eta_k.
$$
Multiplying by $e^{2\nu t}$ and integrating on $[s,s+\bar t]$, we get
$$
\int_{\R^n}
\theta_k(x)|u(\bar t,x)|^2\, d x\leq e^{-2\nu(\bar t-s)}\int_{\R^n}
\theta_k(x)|u(s,x)|^2\, d x+\eta_k{1\over 2\nu}(1-e^{-2\nu (\bar t-s)}),
$$
which in turn implies the thesis.
\end{proof}

Now we can prove
\begin{thm}\label{ascomp1} Assume $(a_{ij}(\tau))_{ij}$ satisfies condition
{\bf(H1)}  and
$F(\tau,x,u)$ satisfies conditions {\bf(H2)}--{\bf(H4)}, {\bf(AP)} and {\bf(D)}. Let $\omega>0$ be fixed,
let $B\subset H^1(\R^N)$ be bounded and let
$(\lambda_n)_{n\in\N}$ be a sequence in $[0,1]$. Then the set $\Sigma\times B$ is
$\{P_{\lambda_n,\omega}\}$-admissible.
\end{thm}
\begin{proof}
First, we chose $R>0$ such that $B\subset B_{H^1}(R;0)$.
By Proposition \ref{localexist}, there exists $T_R>0$
such that, for all $u\in B_{H^1}(R;0)$, for all $\lambda\in[0,1]$, for all $s\in\R$ and for all
$\sigma\in\Sigma$, $\Pi^{\sigma}_{\lambda,\omega}(t,s)u$ is defined for $t\in[s,s+T_R]$ and
$\|\Pi^{\sigma}_{\lambda,\omega}(t,s)u\|_{H^1}\leq 2R$ for $t\in[s,s+T_R]$.
\par Now let $((\sigma_n,u_n))_{n\in\N}$ be a sequence in $\Sigma\times H^1(\R^N)$
and let $(t_n)_{n\in\N}$ be a sequence of positive numbers such that $t_n\to\infty$ as
$n\to\infty$ and $P_{\lambda_n,\omega}(t)(\sigma_n,u_n)\in\Sigma\times B$ for $t\in[0,t_n]$,
$n\in\N$. The latter amounts to saying that $\Pi^{\sigma_n}_{\lambda_n,\omega}(t,0)u_n\in B$
for $t\in[0,t_n]$, $n\in\N$.
\par Since $\Sigma$ is compact, we can assume, without loss of
generality, that there exists $\cl\sigma_\infty\in\Sigma$ such that
$T_\omega(t_n-T_R)\sigma_n\to\cl\sigma_\infty$ and
$T_\omega(t_n)\sigma_n\to T_\omega(T_R)\cl\sigma_\infty=:\sigma_\infty$ as $n\to\infty$.
Moreover, we can assume that there exists $\lambda_\infty\in[0,1]$ such that
$\lambda_n\to\lambda_\infty$ as
$n\to\infty$.
\par Now, since the set
\begin{equation}\label{set}
\{\,\Pi^{\sigma_n}_{\lambda_n,\omega}(t_n-T_R,0)u_n\mid n\in\N\,\}
\end{equation}
is bounded in $H^1(\R^N)$, then passing to a subsequence if necessary, we can
assume that there exists
$\bar u_\infty\in H^1(\R^N)$ such that
$$
\Pi^{\sigma_n}_{\lambda_n,\omega}(t_n-T_R,0)u_n\weakto \bar u_\infty \quad\text{in
$H^1(\R^N)$ as
$n\to\infty$.}
$$
Notice that $\|\bar u_\infty\|_{H^1}\leq R$, so
$\Pi^{\bar\sigma_\infty}_{\lambda_\infty,\omega}(t,0)\bar u_\infty$ is defined for $t\in[0,T_R]$.
We claim
that
$\Pi^{\sigma_n}_{\lambda_n,\omega}(t_n-T_R,0)u_n\to \bar u_\infty $ in the strong
$L^2$-topology.
To this end, it is enough to show that the set (\ref{set})
is relatively compact in the strong $L^2$ topology, or equivalently that it is
totally bounded.

This is a consequence of Lemma \ref{comp1} and of the Rellich Theorem. In fact, for
$n\in\N$ and $k\in\N$ we have
$$
\int_{\R^N}\theta_k(x)\,|(\Pi^{\sigma_n}_{\lambda_n,\omega}(t_n-T_R,0)u_n)(x)|^2\,d x\leq
R^2e^{-2\nu (t_n-T_R)}+\eta_k,
$$
where $\eta_k\to0$ as $k\to\infty$. Let $\epsilon>0$ be fixed. Take $k$ and $n_0$ so large
that $R^2e^{-2\nu (t_n-T_R)}+\eta_k\leq \epsilon$ for all $n\geq n_0$.
Then
\begin{multline}\label{setsum}
\{\,\Pi^{\sigma_n}_{\lambda_n,\omega}(t_n-T_R,0)u_n\mid n\geq n_0\,\}\\
=\{\,\theta_k\,\Pi^{\sigma_n}_{\lambda_n,\omega}(t_n-T_R,0)u_n+
(1-\theta_k)\,\Pi^{\sigma_n}_{\lambda_n,\omega}(t_n-T_R,0)u_n
\mid n\geq n_0\,\}\\
\subset \{\,\theta_k\,\Pi^{\sigma_n}_{\lambda_n,\omega}(t_n-T_R,0)u_n
\mid n\geq n_0\,\}+\{\,(1-\theta_k)\,\Pi^{\sigma_n}_{\lambda_n,\omega}(t_n-T_R,0)u_n
\mid n\geq n_0\,\}\\
\subset B_{L^2}(\epsilon;0)+\{\,(1-\theta_k)\,\Pi^{\sigma_n}_{\lambda_n,\omega}(t_n-T_R,0)u_n
\mid n\geq n_0\,\}.
\end{multline}
The set
$$
\{\,(1-\theta_k)\,\Pi^{\sigma_n}_{\lambda_n,\omega}(t_n-T_R,0)u_n
\mid n\geq n_0\,\}
$$
consists of functions of $H^1(\R^N)$ which are equal to zero outside the ball of radius $\sqrt2k$
in $\R^{N}$.
On the other hand, the $H^1$-norm
of these functions is bounded by a constant depending only on $R$ and $D$. Then, by the Rellich
Theorem, this set is precompact in $L^2(\R^N)$. Hence we
can cover it by a finite number of balls of radius $\epsilon$ in $L^2(\R^N)$. This observation,
together with (\ref{setsum}),
implies that the set (\ref{set}) is totally bounded
and hence precompact in $L^2(\R^N)$. The claim is proved.
\par Finally, by Lemma \ref{contin1}, we have
\begin{multline*}
\Pi^{\sigma_n}_{\lambda_n,\omega}(t_n,0)u_n
=\Pi^{\sigma_n}_{\lambda_n,\omega}(t_n,t_n-T_R)\Pi^{\sigma_n}_{\lambda_n,\omega}(t_n-T_R,0)u_n\\
=\Pi^{T_\omega(t_n-T_R)\sigma_n}_{\lambda_n,\omega}(T_R,0)
\Pi^{\sigma_n}_{\lambda_n,\omega}(t_n-T_R,0)u_n
\to \Pi^{\cl\sigma_\infty}_{\lambda,\omega}(T_R,0)\cl u_\infty\\\text{in $H^1(\R^N)$ as
$n\to\infty$.}
\end{multline*}
Setting $u_\infty:=\Pi^{\cl\sigma_\infty}_{\lambda,\omega}(T_R,0)\cl u_\infty$, it follows that
$\Pi^{\sigma_n}_{\lambda_n,\omega}(t_n,0)u_n\to u_\infty$ in $H^1(\R^N)$ as $n\to\infty$. The
proof is complete.
\end{proof}

\section{Averaging and continuation of invariant sets}
In this section we assume that $(a_{ij}(\tau))_{ij}$ satisfies condition
{\bf(H1)}  and
$F(\tau,x,u)$ satisfies conditions {\bf(H2)}--{\bf(H4)}, {\bf(AP)} and {\bf(D)}.
Let $\lambda_0\in[0,1]$ and $\omega_0>0$ be fixed.   Let
$K_{\lambda_0,\omega_0}\subset\Sigma\times H^1(\R^N)$ be an isolated invariant set of
$P_{\lambda_0,\omega_0}$ and let
$B_{\lambda_0,\omega_0}$ be an isolating neighborhood of $K_{\lambda_0,\omega_0}$. In view of Proposition
\ref{ascomp1}, if  $B_{\lambda_0,\omega_0}$ is bounded, then it is strongly
$P_{\lambda_0,\omega_0}$-admissible. It follows that $K_{\lambda_0,\omega_0}$ is compact (see Theorem I-4.5 in
\cite{ryba}) and its homotopy index
$h(P_{\lambda_0,\omega_0},K_{\lambda_0,\omega_0})$ is well defined.
\par Now we keep $\omega_0$ fixed and we let $\lambda$ run over $[0,1]$. Let
$K_{\lambda,\omega_0}\subset\Sigma\times H^1(\R^N)$ be an isolated invariant set of $P_{\lambda,\omega_0}$
and assume that there exists $B_{\omega_0}\subset\Sigma\times H^1(\R^N)$, such that, for
every $\lambda\in[0,1]$,  $B_{\omega_0}$ is a bounded isolating neighborhood of
$K_{\lambda,\omega_0}$. Then, thanks to Propositions \ref{1contin1} and \ref{ascomp1}, we can apply the
continuation principle I-12.2 of \cite{ryba}. It follows that
$h(P_{\lambda,\omega_0,},K_{\lambda,\omega_0,})$ does not depend on $\lambda$. In particular,
$h(P_{1,\omega_0},K_{1,\omega_0})=h(P_{0,\omega_0},K_{0,\omega_0})$.
\par We have already noticed that
$P_{0,\omega_0}(t)(\sigma,u)=(T_{\omega_0}(t)\sigma,\pi(t)u)$, so
$P_{0,\omega_0}(t)$ is completely decoupled. It follows that, if $K\subset H^1(\R^N)$ is an
isolated invariant set of $\pi(t)$, then $K_{0,\omega_0}:=\Sigma\times K$ is an isolated invariant set
of $P_{0,\omega_0}(t)$. Moreover, by the product formula I-10.6 of \cite{ryba},
\begin{equation}
h(P_{0,\omega_0},K_{0,\omega_0})=h(T_{\omega_0},\Sigma) \land h(\pi,K).
\end{equation}
We recall that, if $(Y,y_0)$ and $(Z,z_0)$ are two pointed spaces, then the {\it smash product}
$(Y,y_0)\land(Z,z_0)$ is the pointed space $(W,w_0)$, where $W:=(Y\times Z)/(Y\times\{z_0\}\cup
\{y_0\}\times Z)$ and $w_0:=[Y\times\{z_0\}\cup
\{y_0\}\times Z]$. In Lemma 1.1 of \cite{ward0} it was proved that if $(Y,y_0)$
is not contractible and
$Z$ is a compact space, then  $(Y,y_0)\land(Z\,\dot\cup\,\{*\},\{*\})$ is not contractible.
\par In the present situation, $\Sigma$ is a compact invariant set of $T_{\omega_0}$
and an isolating neighborhood as well. Actually, $\Sigma$ is an isolating block
with $\Sigma^-=\emptyset$. It follows that
$h(T_{\omega_0},\Sigma)$ is the homotopy type of the pointed space
$(\Sigma\,\dot\cup\,\{*\},\{*\})$. So, if $h(\pi,K)\not=\underline 0$, then
$h(P_{0,\omega_0},K_{0,\omega_0})\not=\underline 0$.
\par Let $K\subset H^1(\R^N)$ be a compact isolated invariant set of $\pi(t)$, with nontrivial Conley
index, and  let $B\subset H^1(\R^N)$ be a bounded isolating neighborhood of $K$. If $\Sigma\times B$ is an
isolating  neighborhood (of $K_{\lambda,\omega_0}$) relative to $P_{\lambda,\omega_0}$ for all
$\lambda\in[0,1]$, then
$$
h(P_{\lambda,\omega_0},K_{\lambda,\omega_0})=h(P_{0,\omega_0},K_{0,\omega_0})
=h(T_{\omega_0},\Sigma) \land h(\pi,K)\not=\underline 0,\quad\lambda\in[0,1].
$$
 In other words,
the isolated invariant set $K$ of $\pi(t)$ can be ``continued" to a family of isolated invariant sets
$K_{\lambda,\omega_0}$ of $P_{\lambda,\omega_0}$, provided one can find a common isolating
neighborhood of the form $\Sigma\times B$, relative to all the $P_{\lambda,\omega_0}$, $\lambda\in[0,1]$.
If the index of
$K$ is nontrivial, the same is true of the index of $K_{\lambda,\omega_0}$.
We stress that, if this is the case, then $K_{\lambda,\omega_0}\not=\emptyset$: this means that
there exist full bounded solutions of (\ref{pbcauchy}) in $B$.
Therefore we are lead to the following question:
\par \emph{given an isolated invariant set $K$ of
$\pi(t)$, is it possible to find a bounded neighborhood
$B$ of $K$ such that $\Sigma\times B$ is an
isolating  neighborhood relative to $P_{\lambda,\omega_0}$ for all
$\lambda\in[0,1]$?}
\par It turns out that the question has a positive answer if $\omega_0$ is sufficiently
large. We need first to prove the following proposition, which ensures a sort of
``singular" admissibility as $\omega\to\infty$.

\begin{prop}\label{ascomp2}
Let $B\subset H^1(\R^N)$ be
a bounded set, let $(\lambda_n)_{n\in\N}$ be a sequence in
$[0,1]$, let $(\sigma_n)_{n\in\N}$ be an arbitrary sequence in
$\Sigma$, let $(\omega_n)_{n\in\N}$ and $(t_n)_{n\in\N}$ be two
sequences of positive numbers, $\omega_n\to\infty$ and
$t_n\to\infty$  as $n\to\infty$, let $(u_n)_{n\in\N}$ be a
sequence in $H^1(\R^N)$ and assume that
$\Pi^{\sigma_n}_{\lambda_n,\omega_n}(t,0)u_n\in B$ for
$t\in[0,t_n]$, $n\in\N$. Then there exists $u_\infty\in H^1(\R^N)$
such that, up to a subsequence,
$$
\Pi^{\sigma_n}_{\lambda_n,\omega_n}(t_n,0)u_n\to u_\infty
$$
in $H^1(\R^N)$ as $n\to\infty$.
\end{prop}
\begin{proof} The proof is
similar to that of Theorem \ref{ascomp1}. First, we chose
$R>0$ such that $B\subset B_{H^1}(R;0)$. By Proposition
\ref{localexist}, there exists $T_R>0$ such that, for all $u\in
B_{H^1}(R;0)$, for all $\lambda\in[0,1]$, for all $\omega>0$, for
all $s\in\R$ and for all $\sigma\in\Sigma$,
$\Pi^{\sigma}_{\lambda,\omega}(t,s)u$ is defined for
$t\in[s,s+T_R]$ and
$\|\Pi^{\sigma}_{\lambda,\omega}(t,s)u\|_{H^1}\leq 2R$ for
$t\in[s,s+T_R]$. Since $B\subset H^1(\R^N)$ is bounded, there
exists $\bar u_\infty\in H^1(\R^N)$ such that, up to a
subsequence,
$$
\Pi^{\sigma_n}_{\lambda_n,\omega_n}(t_n-T_R,0)u_n\weakto \cl
u_\infty\quad\text{in $H^1(\R^N)$ as $n\to\infty$.}
$$
\par Notice that $\|\bar u_\infty\|_{H^1}\leq R$, so $\pi(t)\bar
u_\infty$ is defined for $t\in[0,T_R]$. Like in the proof of
Proposition \ref{ascomp1}, thanks to Lemma \ref{comp1} and to
the Rellich Theorem, we obtain that
$\Pi^{\sigma_n}_{\lambda_n,\omega_n}(t_n-T_R,0)u_n\to \bar
u_\infty $ in the strong $L^2$-topology. Finally, by Theorem
\ref{contin2}, we have
\begin{multline*}
\Pi^{\sigma_n}_{\lambda_n,\omega_n}(t_n,0)u_n
=\Pi^{\sigma_n}_{\lambda_n,\omega_n}(t_n,t_n-T_R)\Pi^{\sigma_n}_{\lambda_n,\omega_n}(t_n-T_R,0)u_n\\
=\Pi^{T_{\omega_n}(t_n-T_R)\sigma_n}_{\lambda_n,\omega_n}(T_R,0)
\Pi^{\sigma_n}_{\lambda_n,\omega_n}(t_n-T_R,0)u_n \to \pi(T_R)\cl
u_\infty\\\text{in $H^1(\R^N)$ as $n\to\infty$.}
\end{multline*}
Setting $u_\infty:=\pi(T_R)\cl u_\infty$, it follows that
$\Pi^{\sigma_n}_{\lambda_n,\omega_n}(t_n,0)u_n\to u_\infty$ in
$H^1(\R^N)$ as $n\to\infty$. The proof is complete.
\end{proof}

We recall the following

\begin{defn}\label{fullbounded} A curve $t\mapsto u(t)\in
H^1(\R^N)$, $t\in\R$ is said to be a full solution of the process
$\Pi^\sigma_{\lambda,\omega}(t,s)$ iff
$$
u(t)=\Pi^\sigma_{\lambda,\omega}(t,s)u(s)
\quad\text{for all $t\geq s$, $s\in\R$.}
$$
\end{defn}

Now we have:

\begin{cor}\label{fullbound} Let $B\subset H^1(\R^N)$,
$(\lambda_n)_{n\in\N}$, $(\sigma_n)_{n\in\N}$ and
$(\omega_n)_{n\in\N}$ be as in Proposition \ref{ascomp2}. For all
$n\in\N$, let $u_n\colon\R\to H^1(\R^n)$ be a full solution of
$\Pi^{\sigma_n}_{\lambda_n,\omega_n}(t,s)$, such that $u_n(t)\in B$
for all $t\in\R$. Under these hypotheses, there exists a
subsequence of $(u_n)_{n\in\N}$, again denoted by
$(u_n)_{n\in\N}$, and a full solution $u_\infty\colon\R\to
H^1(\R^n)$ of the averaged semiflow $\pi(t)$, such that
$u_n(t)\to u_\infty (t)$ in $H^1(\R^N)$ as $n\to\infty$, uniformly on
every bounded subinterval of $\R$.
\end{cor}
\begin{proof} As in the proof of Proposition
\ref{ascomp2}, we begin by taking $R>0$ such that $B\subset
B_{H^1}(R;0)$. By Proposition \ref{localexist}, there exists
$T_R>0$ such that, for all $u\in B_{H^1}(R;0)$, for all
$\lambda\in[0,1]$, for all $\omega>0$, for all $s\in\R$ and for
all $\sigma\in\Sigma$, $\Pi^{\sigma}_{\lambda,\omega}(t,s)u$ is
defined for $t\in[s,s+T_R]$ and
$\|\Pi^{\sigma}_{\lambda,\omega}(t,s)u\|_{H^1}\leq 2R$ for
$t\in[s,s+T_R]$. Next, we fix once and for all a sequence
$(t_n)_{n\in\N}$ of positive numbers, with $t_n\to\infty$ as
$n\to\infty$.
\par Let $k\in\Z$. For all sufficiently large $n$, we have
\begin{multline*}
u_n(kT_R)=\Pi^{\sigma_n}_{\lambda_n,\omega_n}(kT_R,kT_R-t_n)u_n(kT_R-t_n)\\
=\Pi^{T_{\omega_n}(kT_R-t_n)\sigma_n}_{\lambda_n,\omega_n}(t_n,0)u_n(kT_R-t_n).
\end{multline*}
Then, by Theorem \ref{ascomp2}, there is a subsequence of
$(u_n(kT_R))_{n\in\N}$, again denoted by $(u_n(kT_R))_{n\in\N}$,
and there exists $v(kT_R)\in H^1(\R^n)$ such that $u_n(kT_R)$
converges strongly to $v(kT_R)$ in $H^1(\R^n)$ as $n\to\infty$. In
particular, $\|v(kT_R)\|_{H^1}\leq R$. Using Cantor's diagonal
procedure we obtain the existence of a subsequence of
$(u_n)_{n\in\N}$, again denoted by $(u_n)_{n\in\N}$,
and a sequence $v(kT_R)\in H^1(\R^n)$, $k\in\Z$, such that, for
every $k\in \Z$,
$$
u_n(kT_R)\to v(kT_R)\quad\text{in $H^1(\R^n)$ as $n\to\infty$.}
$$
By Theorem \ref{contin2}, we have that, for
all $k\in\Z$,
$$
\Pi^{\sigma_n}_{\lambda_n,\omega_n}(t,kT_R)u_n(kT_R) \to\pi(t-kT_R)v(kT_R)
$$
in $H^1(\R^n)$ as $n\to\infty$, uniformly on $[kT_R,(k+1)T_R]$.
\par In particular, one has
$\Pi^{\sigma_n}_{\lambda_n,\omega_n}((k+1)T_R,kT_R)u_n(kT_R)\to\pi(T_R)v(kT_R)$.
On the other hand,
$\Pi^{\sigma_n}_{\lambda_n,\omega_n}((k+1)T_R,kT_R)u_n(kT_R)=u_n((k+1)T_R)\to
v((k+1)T_R)$. Hence we deduce that $v((k+1)T_R)=\pi(T_R) v(kT_R)$
for all $k\in\Z$. We can therefore define
$$
u_\infty(t):=\pi(t-kT_R)v(kT_R)\quad \text{for $t\in [kT_R,(k+1)T_R]$,}
$$
which is easily seen to be a full solution of of $\pi(t)$.
Moreover,
$$
u_n(t)\to u_\infty(t)\quad\text{as $n\to\infty$}
$$
uniformly on every bounded subinterval of $\R$. \end{proof}

Finally, we can prove:
\begin{thm}\label{singular}
Let $K$ be an isolated invariant set of $\pi(t)$ and let $B\subset H^1(\R^N)$
be a bounded isolating neighborhood of $K$. There exists $\bar\omega>0$ such that, for all
$\omega>\bar\omega$ and for all $\lambda\in[0,1]$, $\Sigma\times B$ is an isolating neighborhood
relative to $P_{\lambda,\omega}$.
\end{thm}
\begin{proof} Assume by contradiction that the theorem is not true. Then there exist
a sequence $(\lambda_n)_{n\in\N}$ in $[0,1]$, a sequence of positive numbers $(\omega_n)_{n\in\N}$,
$\omega_n\to +\infty$ as $n\to\infty$, a sequence
$(\sigma_n)_{n\in\N}$ in $\Sigma$ and
a sequence $(u_n)_{n\in\N}$ of functions from $\R$ to $H^1(\R^N)$, such that, for
$n\in\N$,  $u_n(t)$ is a full solution of
$\Pi^{\sigma_n}_{\lambda_n,\omega_n}(t,s)$, with $u_n(t)\in B$
for all $t\in\R$ and $u_n(0)\in\partial B$ for all $n\in\N$.
By Corollary \ref{fullbound},
there exists a
subsequence of $(u_n)_{n\in\N}$, again denoted by
$(u_n)_{n\in\N}$, and a full solution $u_\infty\colon\R\to
H^1(\R^n)$ of the averaged semiflow $\pi(t)$, such that
$u_n(t)\to u_\infty (t)$ as $n\to\infty$ uniformly on
every bounded subinterval of $\R$. It follows that $u_\infty(t)\in B$
for all $t\in\R$ and $u_\infty(0)\in\partial B$, thus contradicting the fact that $B$ is an isolating
neighborhood relative to $\pi(t)$.
\end{proof}

The results proved in this section can be summarized as follows:

\begin{thm}\label{averconley} Assume that $(a_{ij}(\tau))_{ij}$ satisfies condition
{\bf(H1)}  and
$F(\tau,x,u)$ satisfies conditions {\bf(H2)}--{\bf(H4)}, {\bf(AP)} and {\bf(D)}.
Suppose that the semiflow $\pi(t)$, generated by the
autonomous averaged equation (\ref{pbcauchyaut}), posseses an
isolated invariant set $K\subset H^1(\R^N)$, with nontrivial
homotopy index. Then, for all sufficiently large $\omega$ and for all
$\lambda\in[0,1]$,
the skew-product semiflow generated by the non-autonomous equation
(\ref{pbcauchy}) possesses an isolated invariant set $K_{\lambda,\omega}\subset
\Sigma\times H^1(\R^N)$, with nontrivial homotopy index.\end{thm}

\section{Recurrent motions}
In this section we shall discuss some consequences of Theorem
\ref{averconley}. Let $\lambda=1$. If
$h(P_{1,\omega},K_{1,\omega})\not=\underline 0$, then
$K_{1,\omega}\not=\emptyset$. This means that there exist
$(\sigma_0, u_0)\in \Sigma\times H^1(\R^N)$ and a function
$(\sigma,u)\colon\R\to\Sigma\times H^1(\R^N)$, such that
$(\sigma(0),u(0))=(\sigma_0,u_0)$, $(\sigma(t),u(t))\in
K_{1,\omega}$ for all $t\in\R$ and
$(\sigma(t),u(t))=P_{1,\omega}(t-s)(\sigma(s),u(s))$ for all
$t\geq s$. It follows that $u(t)$ is a bounded full solution of
the process $\Pi^{\sigma_0}_{1,\omega}$. If we are interested in
proving the existence of bounded full solutions of the original
equation (\ref{equazione}), we can argue as follows. Since the
orbit $\{\sigma(t)\mid t\in\R\}$ is dense in $\Sigma$, then there
exists a sequence $(t_n)_{n\in\N}$, such that
$\sigma(t_n)\to\sigma_\sharp:=((a_{ij})_{ij}, F)$ as $n\to\infty$.
Since $K_{1,\omega}$ is compact, we can assume, without loss of
generality, that there exists $u_\sharp\in H^1(\R^N)$ such that
$(\sigma_\sharp,u_\sharp)\in K_{1,\omega}$ and $u(t_n)\to
u_\sharp$ as $n\to\infty$. It follows that there exists a function
$(\tilde\sigma,\tilde u)\colon\R\to\Sigma\times H^1(\R^N)$, such
that $(\tilde\sigma(0),\tilde u(0))=(\sigma_\sharp,u_\sharp)$,
$(\tilde\sigma(t),\tilde u(t))\in K_{1,\omega}$ for all $t\in\R$
and $(\tilde\sigma(t),\tilde
u(t))=P_{1,\omega}(t-s)(\tilde\sigma(s),\tilde u(s))$ for all
$t\geq s$. It follows that $\tilde u(t)$ is a bounded full
solution of the process $\Pi^{\sigma_\sharp}_{1,\omega}$, i.e. a
bounded full solution of (\ref{equazione}).
\par From the dynamical point of view, it is much more interesting to look for {\it recurrent
solutions} rather than for {\it bounded solutions} of the equation (\ref{equazione}).
\par Let $X$ be a complete metric space and let $\pi(t)$ be a {\it global two-sided flow} on $X$. The
following basic concepts were introduced by Birkhoff (see \cite{birk}; for a modern
treatment, see also the book of Bhatia and Szeg\"o \cite{BS}):

\begin{defn}\label{recurr} A point $x\in X$ is called {\rm recurrent} iff
\begin{enumerate}
\item the orbit $\{\pi(t)x\mid t\in\R\}$ is precompact in $X$;
\item for every $\epsilon>0$ there exists $\ell>0$ such that in every interval $I\subset \R$ of lenght
$\ell$ there is a $\tau$ such that $d(\pi(\tau)x,x)<\epsilon$.
\end{enumerate}
If the point $x$ is recurrent, the same is true of the point $\pi(t)x$, for all $t\in\R$. The full
trajectory $\pi(\cdot)x$ is then called {\rm recurrent}.
\end{defn}

\begin{defn}\label{minimal}
A set $M\subset X$ is called a {\rm minimal set} iff
\begin{enumerate}
\item $M$ is closed and invariant;
\item $M$ does not contain nonempty, proper, closed invariant subsets.
\end{enumerate}
\end{defn}

The concepts of {\it recurrent point} and {\it minimal set} are related by the following theorem (for a
proof, see e.g. \cite{BS}):

\begin{thm}[Birkhoff, 1926]\label{thbirk1} A point $x\in X$ is recurrent if and only
if it belongs to a compact minimal set.\end{thm}

The existence of recurrent points for a flow in a compact metric space is guaranteed by the following

\begin{thm}[Birkhoff, 1926]\label{thbirk2} If $X$ is compact, then there exists a
minimal set $M\subset X$.\end{thm}

Concerning the semiflow $P_{1,\omega}$, we stress that its phase space is not compact. Moreover,
the trajectories are in general defined only in forward time. However, we can restrict the semiflow
$P_{1,\omega}$ to the compact invariant set $K_{1,\omega}$. Notice that, for every $(\sigma_0,u_0)\in
K_{1,\lambda}$, there is a function $(\sigma,u)\colon\R\to\Sigma\times H^1(\R^N)$,
such that $(\sigma(0),u(0))=(\sigma_0,u_0)$, $(\sigma(t),u(t))\in K_{1,\omega}$ for all $t\in\R$ and
$(\sigma(t),u(t))=P_{1,\omega}(t-s)(\sigma(s),u(s))$ for all $t\geq s$. If the semiflow
$P_{1,\omega}$, restricted to $K_{1,\omega}$, possesses the {\it backward uniqueness} property, then
it admits a unique {\it flow extension}. Thanks to an abstract result of Lions and Malgrange (\cite{liomal}),
the {\it backward uniqueness} property holds for equation (\ref{pbcauchy}), provided
we replace the H\"older condition (\ref{a0}) for $a_{ij}(\cdot)$  in {\bf (H1)} with the following stronger
Lipschitz condition: for all $\tau_1,\tau_2\in\R$, and for $1\leq i,j\leq N$,
\begin{equation}\label{lip}
|a_{ij}(\tau_1)-a_{ij}(\tau_2)|\leq C|\tau_1-\tau_2|.
\end{equation}

Under this stronger assumption, we can apply Birkhoff's theorem to the unique {\it flow extension}
of the semiflow $P_{1,\omega}$ in the compact metric space $K_{1,\omega}$. We thus obtain
the existence of a {\it minimal set} $M_{1,\omega}$ contained in $K_{1,\omega}$. This in turn implies the
existence  of at least one recurrent trajectory in $K_{1,\omega}$.
\par To the concept of {\it recurrent trajectory} there corresponds the concept of {\it recurrent
function}. Let $Y$ be a complete metric space and let ${\mathcal U}(\R,Y)$ be the space of all continuous
functions $g\colon\R\to Y$, with the (metrizable) topology of uniform convergence on the bounded segments.
For $g\in {\mathcal U}(\R,Y)$ and $s\in\R$, define $(T(s)g)(t):=g(t+s)$, $t\in\R$.
A function $g\in {\mathcal U}(\R,Y)$
is called {\it recurrent} if the trajectory $T(s)g$ is recurrent in ${\mathcal U}(\R,Y)$.
\par The connection between
recurrent functions and recurrent trajectories is the following: let $x$ be a recurrent point of a global
flow $\pi(t)$ in a complete metric space $X$;  let $Y$ be a complete metric space and let $\phi\colon X\to
Y$ be a continuous function; then the function $t\mapsto \phi(\pi(t)x)$ is recurrent.
Therefore, if $(\sigma_0, u_0)$ is a recurrent point of the flow extension of
$P_{1,\omega}$ in  $K_{1,\omega}$, then there exists a recurrent function
$(\sigma,u)\colon\R\to\Sigma\times H^1(\R^N)$, such that $(\sigma(0),u(0))=(\sigma_0,u_0)$,
$(\sigma(t),u(t))\in K_{1,\omega}$ for all
$t\in\R$ and $(\sigma(t),u(t))=P_{1,\omega}(t-s)(\sigma(s),u(s))$ for all $t\geq s$.
It follows that $u(t)$ is recurrent solution of the process
$\Pi^{\sigma_0}_{1,\omega}$.

If we are interested in proving the existence of recurrent solutions of the
original equation (\ref{equazione}), we can argue as follows. Since the orbit
$\{\sigma(t)\mid t\in\R\}$ is dense in
$\Sigma$, then there exists a sequence $(t_n)_{n\in\N}$, such that
$\sigma(t_n)\to\sigma_\sharp:=((a_{ij})_{ij}, F)$ as $n\to\infty$. Since $M_{1,\omega}$ is compact, we can
assume, without loss of generality, that there exists $u_\sharp\in H^1(\R^N)$ such that
$(\sigma_\sharp,u_\sharp)\in M_{1,\omega}$ and $u(t_n)\to u_\sharp$ as $n\to\infty$. It follows that here
exists a function
$(\tilde\sigma,\tilde u)\colon\R\to\Sigma\times H^1(\R^N)$, such that
$(\tilde\sigma(0),\tilde u(0))=(\sigma_\sharp,u_\sharp)$,
$(\tilde\sigma(t),\tilde u(t))\in M_{1,\omega}$ for all
$t\in\R$ and
$(\tilde\sigma(t),\tilde u(t))=P_{1,\omega}(t-s)(\tilde\sigma(s),\tilde u(s))$ for all $t\geq s$. It
follows that
$\tilde u(t)$ is a recurrent solution of the process $\Pi^{\sigma_\sharp}_{1,\omega}$, i.e. a
recurrent solution of \ref{equazione}. We can summarize the above considerations in the following

\begin{thm}\label{recurrsol} Assume that $(a_{ij}(\tau))_{ij}$ satisfies condition
{\bf(H1)}, with the H\"older condition (\ref{a0}) replaced by the Lipschitz condition (\ref{lip}),  and
$F(\tau,x,u)$ satisfies conditions {\bf(H2)}--{\bf(H4)}, {\bf(AP)} and {\bf(D)}.
Suppose that the semiflow $\pi(t)$, generated by the
autonomous averaged equation (\ref{pbcauchyaut}), possesses an
isolated invariant set $K\subset H^1(\R^N)$, with nontrivial
homotopy index. Then, for all sufficiently large $\omega$, the non-autonomous equation
(\ref{equazione}) possesses a recurrent solution.\end{thm}

\par We conclude with an example, in which the averaged equation is asymptotically linear (cf \cite{priz}).
More precisely, we assume that the average $\bar F(x,u)$ satisfies (\ref{a1}) with $\beta=0$
and (\ref{a4}) with $q=2$. Moreover, we assume that
\begin{equation}\label{aslin}
\lim_{|u|\to\infty}{{\bar F(x,u)}\over{u}}=V(x):=-V_1(x)+V_2(x)\quad
\text{for all $x\in\R^n$,}\end{equation}
where $V_1\in L^\infty(\R^n)$, with
$V_1(x)\geq\tilde\nu>0$ for all $x\in\R^n$, and $V_2\in L^{\rho}(\R^n)$,
with $n\leq \rho<\infty$.
It was observed in \cite{priz} that the essential spectrum of the
operator
$-\Delta-V(\cdot)$ is contained in $[\tilde\nu,+\infty[$. In
particular, the part of the spectrum of $-\Delta-V(\cdot)$
contained in $]-\infty, \tilde\nu/2[$ is a finite set, consisting
of isolated eigenvalues with finite multiplicity.
We assume that the
following {\it non-resonance condition at infinity} is  satisfied:
\begin{equation}\label{nrc}
{\rm{ker}}(-\Delta-V(\cdot))=(0).
\end{equation}
In \cite{priz} it was proved the following

\begin{thm}\label{Sigma} Assume that $\bar F$ satisfies (\ref{a1}) with $\beta=0$, (\ref{a4}) with $q=2$,
(\ref{aslin}) and (\ref{nrc}). Let
$m$ be the total multiplicity of the negative eigenvalues of
$-\Delta-V(\cdot)$. Denote by $\pi_{\bar F}$ the semiflow generated by (\ref{pbcauchyaut}) and by
$K_{\bar F}$ the union of the ranges of all bounded full solutions of $\pi_{\bar F}$. Then $K_{\bar F}$
is a compact isolated invariant set with homotopy index
$$
h(\pi_F,K_{\bar F})=\Sigma^m,
$$
where $\Sigma^m$ is the homotopy type of a $m$-dimensional pointed sphere. In particular,
$h(\pi_F,K_{\bar F})\not=\underline 0$, so $K_{\bar F}\not=\emptyset$.
\end{thm}

From Theorems \ref{recurrsol} and \ref{Sigma} one can finally
deduce:

\begin{thm}\label{recurrsol2} Assume that $(a_{ij}(\tau))_{ij}$ satisfies condition
{\bf(H1)}, with the H\"older condition (\ref{a0}) replaced by the Lipschitz condition (\ref{lip}),  and
$F(\tau,x,u)$ satisfies conditions {\bf(H2)}--{\bf(H4)}, {\bf(AP)} and {\bf(D)}.
Assume that the average $\bar F$ satisfies (\ref{a1}) with $\beta=0$, (\ref{a4}) with $q=2$,
(\ref{aslin}) and (\ref{nrc}). Then, for all sufficiently large $\omega$, the non-autonomous equation
(\ref{equazione}) possesses a recurrent solution.\end{thm}

\smallskip


\end{document}